\documentclass[12pt]{article}

\input epsf

 \makeatletter
\renewcommand\section{\@startsection{section}{1}{\z@}%
       {-2.5ex \@plus -1ex \@minus -.2ex}%
        {1.5ex \@plus.2ex}%
        {\reset@font\large\bfseries}}
\renewcommand\subsection{\@startsection{subsection}{2}{\z@}%
                                     {-1.5ex\@plus -0.5ex \@minus -.2ex}%
                                     {1ex \@plus .2ex}%
                                     {\reset@font\normalsize\bfseries}}

\renewcommand\maketitle{
  \vskip 30\p@
  \begin{center}%
    {\large \bf\@title \par}%
    \vskip 1em%
    {\large
     \lineskip .75em%
      \begin{tabular}[t]{c}%
        \small \@author 
      \end{tabular}\footnote{}
     \par}}
  \global\let\thanks\relax
  \global\let\maketitle\relax
  \global\let\@thanks\@empty
  \global\let\@author\@empty
  \global\let\@date\@empty
  \global\let\@title\@empty
  \global\let\title\relax
  \global\let\author\relax
  \global\let\date\relax
  \global\let\and\relax

\makeatother

\begin{document}

\def \tvi {\vrule height 12pt depth 5pt width 0pt}

\def \pc #1#2|{{\uppercase{#1}}%
                {\sc {#2}}}
\def \pd #1 {\pc #1| \space}

\def \oe {o\hskip -1.8pt e}

\maketitle 

\setcounter{page}{1}

\thispagestyle{plain}
 \medskip

\font\petcap=cmcsc10
\def\pointir{\unskip . --- \ignorespaces}
\def\noipc#1{\noindent{\petcap #1}}
\def\zed {\hbox {\rm { Z \kern -2.8ex Z}\kern 1.15ex}}
\def\real {\hbox  {\rm {R \kern -2.8ex I}\kern 1.15ex}}
\def\integer {\hbox {\rm {N \kern -2.8ex I}\kern 1.15ex}}
\def\hyp {\hbox  {\rm {H \kern -2.8ex I}\kern 1.15ex}}
\def\cqfd{\penalty 500 \hfill\hfill \hbox{ \vrule height6pt width4pt depth1pt}}
\def\dem{\noipc{Proof. }}
\def\dist{\hbox{dist}}
\def\Isom{\hbox{Isom}}
\def\interieur#1{\mathord{\mathop{\kern 0pt #1}\limits^\circ}}
\def\demi{{\ts{1\over 2}}}
\def\chapeau#1{\smash{\widehat#1}}
\def\norme#1{{\|#1\|}}\def\arret#1{^{|{#1}}}
~~\vskip .7in

\centerline {\bf The  Weil-Petersson K\"ahler form and affine foliations
on surfaces } 

\centerline {\bf }

\vskip .1in
\centerline{by}
\medskip
\centerline {Athanase Papadopoulos and R. C. Penner}

\vskip .3in

\leftskip .5in\rightskip .5in
\noindent ${\underline{\rm Abstract}}$~~The space of broken hyperbolic structures generalizes
the usual Teichm\"uller space of a punctured surface, and the space of projectivized broken measured foliations--or equivalently, the space
of projectivized affine foliations of a punctured surface--likewise admits a generalization to projectivized broken measured foliations.
Just as projectivized measured foliations provide Thurston's boundary for Teichm\"uller space, so too do projectivized broken
measured foliations provide a boundary for the space of broken hyperbolic structures.  In this paper, we naturally extend the Weil-Petersson
K\"ahler two-form to a corresponding two-form on the space of broken hyperbolic structures as well as Thurston's symplectic form to a
corresponding two-form on the space of broken measured foliations, and we show that the former limits in an appropriate sense to the latter. 
The proof in sketch follows earlier work of the authors for measured foliations and depends upon techniques from decorated
Teichm\"uller theory, which is also applied here to a further study of broken hyperbolic structures. 

\leftskip=0ex\rightskip=0ex
\vskip .3in

\noindent  $\underline{ \bf
\hbox{0.--- Introduction}}$ 
\vskip .2in

\medskip
In the paper [PP], we established a relation between the    
Weil-Petersson K\"ahler form
 on the Teichm\"uller space of a punctured surface and
Thurston's piecewise-linear symplectic form on the space of measured foliations
of compact support on that surface. 
We extend here this work to  
the context of broken hyperbolic structures
and of broken measured foliations.  Both spaces of broken structures
were  defined in [OP2], and  we shall recall the definitions in the next
section. The space of broken hyperbolic structures contains Teichm\"uller
space as a proper subspace; the space 
of 
 broken
measured foliations can be identified with the space
of affine foliations on the surface (as discussed in Appendix~A), and it
likewise contains  the space of measured foliations as a proper subspace. 

More precisely, in this paper, we define a two-form 
on the space of broken hyperbolic structures
which
extends   the Weil-Petersson K\"ahler two-form defined on  Teichm\"uller
space.
Likewise, we define a two-form
on the space of broken
measured foliations, which extends Thurston's form 
defined on measured foliations space.  In the paper [OP2],   the
space of projective classes of broken
measured foliations was described as  a boundary for
the space of broken hyperbolic structures, generalizing 
Thurston's realization of
  the space of projective classes of measured
foliations  as a boundary to  the Teichm\"uller space of the surface.
In the work here, we exhibit  a relation between the  two-forms on the spaces of
broken hyperbolic structures and broken measured foliations, which is
analogous to the
relation on Teichm\"uller space and its boundary that we produced 
 in [PP].  In fact,
the whole discussion in this paper 
takes place (as in the paper [PP]) in the {\it decorated}
spaces of broken hyperbolic structures and of broken measured foliation,
whose definitions we give in the next section.  
Finally in Appendix~B, we further apply techniques from decorated Teichm\"uller theory to study the space of
broken hyperbolic structures.

\vskip .3in

\noindent  $\underline{ \bf
\hbox{1.---  Decorated broken hyperbolic 
structures and broken measured foliations}}$ 
\vskip .2in

In what follows, $F_g^s$ is an oriented surface of negative Euler characteristic
with genus $g$ and $s$ punctures,
$s>0$.  
Let $\Delta$ be an ideal triangulation of $F_g^s$, that is, a
decomposition of $F_g^s$ into triangles whose  vertices are  at the 
punctures.
 (We do not require
that two triangles meet in at most one edge, as in   a
simplicial
decomposition.) Let $\widetilde{F_g^s}$ be the topological universal cover of $F_g^s$ 
(with respect to some basepoint which we suppress) and  
$\widetilde{\Delta}$ the lift of $\Delta$ to $\widetilde{F_g^s}$.

\medskip

\noindent {\bf Definition 1.1 (Broken hyperbolic structures).---} 
 A {\it broken hyperbolic structure} on $(F_g^s,\Delta)$
  is a  
 Riemannian metric of constant curvature $-1$  on
 $F_g^s-\Delta$   such that the completion of each face of
 $\Delta$, with its vertices deleted, is isometric 
    to    a hyperbolic  ideal triangle.    
  (We recall that a hyperbolic  ideal triangle
   is the convex hull, in the hyperbolic plane $\hyp^2$, of three distinct
   points in the boundary of that space.) We require furthermore that   
  conditions (1.1.1) to (1.1.3) be satisfied :
   
   \medskip
   
   \noindent (1.1.1) {(\bf Homothety)} \hskip .1in  Consider the surface
$F_g^s$
   as a quotient of the   
disjoint union of the ideal triangles
  (the faces of $\Delta$),
by  
 gluing their edges pairwise. 
Then, each gluing map between two edges  
  of triangles which map  to
the same edge   of $\Delta$ is a homothety, with respect to    the metrics on the
edges of the triangles which are induced by the  hyperbolic metrics on the triangles; that is, across each pair of identified edges, the
hyperbolic metrics on the faces are related by an overall homothety.  Thus, associated with a pair $f_1,f_2$ of faces sharing the edge $e$,
there is a ``homothetic scaling factor'' $\sigma (f_1,f_2)\in\real_+$ of broken metric across
$e$; specifically, if an arc $a\subseteq e$ has length in $f_i$ given by $\rho _i$, for $i=1,2$,
then $\rho _2=\sigma (f_1,f_2)~\rho_1$.  Notice that $\sigma (f_1,f_2)~\sigma (f_2,f_1)=1$.

 \medskip

 \noindent (1.1.2) {\bf (No Horocyclic Holonomy)}\hskip .1in  {The
holonomy of the broken  hyperbolic structure 
  is   trivial around each cusp}. This means that
if $f_1,f_2,\ldots , f_n,f_{n+1}=f_1$ are consecutive
faces of $\Delta$ traversed in a counterclockwise sense about a puncture, then 
$\prod _{i=1}^n
\sigma (f_i,f_{i+1})=1$.  (In particular, for $s=1$, this condition holds
automatically.)

\medskip

   \noindent (1.1.3) {\bf (Completeness)}\hskip .1in  The broken 
hyperbolic 
   structures that we consider in this paper
   are {complete} in the sense which we describe now. To formulate the
   definition, we first 
     introduce
   the
   {\it horocyclic foliation} associated to a broken hyperbolic structure.
   This is a   foliation on $F_g^s$ which is defined as follows.
   On each  ideal triangle in $\hyp^2$, there is a well-defined
     partial foliation  (that is, a foliation supported on a subsurface
	of that triangle)
     whose leaves are connected pieces of horocycles or ``horocyclic
arcs'', where the horocycle is 
     centered at the vertices of these triangles and where each leaf has its  endpoints
        on
    the edges of that triangle making  right angles
    with the  edges. The
    non-foliated region is a central  
    triangle which is bounded by three of these horocyclic arcs, which
pairwise 
    meet tangentially
    at their endpoints (Figure 1).

\medskip


~~\vskip 2.2in

\hskip .9in\epsffile{WP1ai.epsf}

\vskip .1in

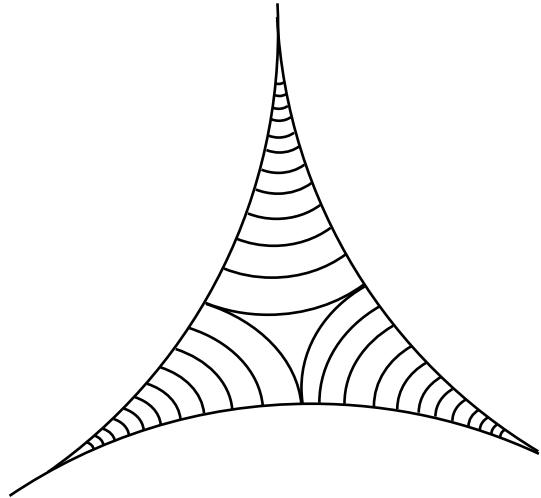
\centerline{{\bf Figure 1} The horocyclic partial foliation of an ideal
triangle.}

\medskip

We  call this foliation the  
    {\it
    horocyclic (partial) foliation  of the ideal triangle}.
    (We note that since all ideal triangles are
   isometric, it suffices  to define   the horocyclic 
    foliation of a particular 
   ideal triangle, and then carry by an isometry to an arbitrary triangle.)
      From this partial
     foliation we can obtain a   foliation of full support 
     on the ideal triangle, by collapsing the
     non-foliated region onto a {\it tripod} (see Figure 2).  

\medskip


~~\vskip 2.2in

\hskip .9in\epsffile{WP2ai.epsf}

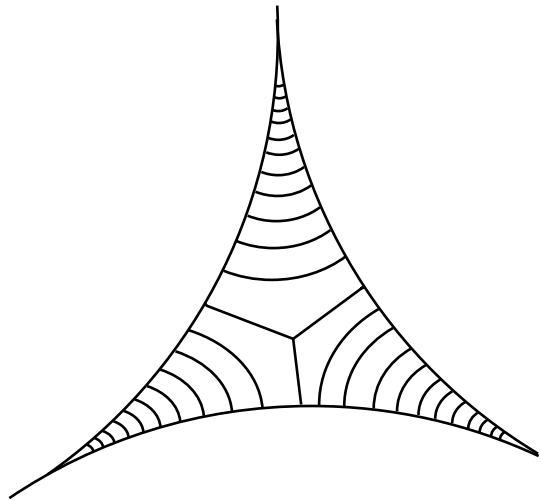
\centerline{{\bf Figure 2} The horocyclic foliation of an ideal
triangle.}

\medskip

The
    horocyclic foliations of the various ideal triangles
    which constitute  the faces of $\Delta$ extend  naturally
   to  a  foliation on $F_g^s$ which we call the {\it horocyclic
   foliation  of the broken hyperbolic structure}.
   By construction, this
     foliation  is
     well-defined up to an isotopy preserving the edges of $\Delta$.  
    The broken hyperbolic structure
   is then said to be  {complete} if there exists a
    neighborhood of each puncture of
   $F_g^s$   which is topologically 
   an annulus and on which the  foliation 
  induced by the horocyclic foliation  is a foliation by  circles
   which are homotopic to the puncture.

\medskip

   Let us note that in the case where the broken hyperbolic structure on 
   $(F_g^s,\Delta)$ is
   a {\it hyperbolic  structure} (in other words, if all the homothety scales in (1.1.1)
   are equal to unity), then completeness of 
   broken hyperbolic structure  is equivalent  to completeness of hyperbolic
   structure in the usual sense. We recall that 
   in this case, the neighborhood  of each puncture is a
   ``cusp", that is, a   surface isometric to the quotient of the region
   in the upper half space model of $\hyp^2$ which is above the line $y=1$ by
    a parabolic
   transformation fixing the point $\infty$.

   \medskip

\noindent {\bf Equivalence relation.---}   
  Two broken hyperbolic
structures  on
$(F_g^s,\Delta)$ are considered to be  equivalent if they
differ by an isometry which is isotopic to the identity, with this 
isotopy   preserving   setwise the vertices and the edges of  $\Delta$. 

\medskip

Let ${\cal BH}(\Delta)$ denote the set of equivalence classes of broken
hyperbolic structures.

\medskip

The subset of ${\cal BH}(\Delta)$ which consists of the broken hyperbolic
structures
 for
which the homothety scales in (1.1.1) are all equal to unity can naturally be identified with the
Teichm\"uller space 
${\cal T}= {\cal T}(F_g^s)$ of   $F_g^s$, 
that is, the space of isotopy classes of complete finite-area hyperbolic
metrics on this surface. Thus, we have a natural inclusion ${\cal T}\subset
{\cal BH}(\Delta)$.

 \medskip

\noindent {\bf Definition 1.2 (Decorated broken hyperbolic structure).---} 
  A
{\it decorated broken hyperbolic structure} on $(F_g^s,\Delta)$ 
is a broken hyperbolic structure
together with the choice, for each puncture of $F_g^s$,  of a closed
leaf of the associated horocyclic  foliation
 which is homotopic to that puncture 
 (that is, a closed leaf contained in one of the annuli
 that are discussed in the construction of horocyclic foliation following
Property 1.1.2  of Definition 1.1).

 \medskip

The space of decorated broken hyperbolic structures, 
  up to  the equivalence relation generated by  isometries   isotopic to the identity and
fixing setwise the vertices and the edges of  $\Delta$,
 will be denoted by
  $\widetilde{{\cal BH}}(\Delta)$.

\medskip

 The decoration of the surface (that is, the choice of a horocycle around each
 puncture) is a tool which has been proved to be useful in hyperbolic
 geometry (see [Pe1]). In some sense, it serves
 as a way for measuring (algebraic) distances from the
 puncture using the Minkowski inner product, cf. Appendix B.

\medskip

The {\it decorated Teichm\"uller space} $\widetilde{\cal T}$ from [Pe1] is 
the total space of the trivial $\real_+^s$ bundle over 
${{\cal T}}$, where the fibre over a point is the collection of all $s$-tuples of horocycles,
not necessarily embedded or disjoint, with one horocycle about each puncture; the coordinate on the fibre
is the hyperbolic length of the specified horocycle.  In this paper, we shall consider
the subspace $\widetilde{\cal T}'\subset\widetilde{\cal T}$ corresponding to {\sl disjointly embedded} families of horocycles.
Thus, whereas $\widetilde{\cal T}\not\subset\widetilde{BH}(\Delta )$, we have the natural inclusion $\widetilde{\cal
T}'\subset\widetilde{BH}(\Delta )$.

 \medskip
 
\noindent  {\bf Definition 1.3 (Measured foliation).---} A foliation $F$ on  a
surface
is said to be a {\it measured foliation} if   each arc  $c$ which is transverse
to the leaves of $F$ is equipped with a Borel measure which is equivalent to
the Lebesgue measure of an interval, with  the property
 that if   $c'$ is another arc which is  transverse to $F$ and which is
  obtained from $c$ by an isotopy
during  which each point of the arc stays on the same leaf, then the resulting natural map between
$c$ and $c'$ (obtained by sliding along the leaves)  is measure-preserving.

\medskip

\noindent {\bf Definition 1.4. (Broken measured foliation).---}  A {\it broken
measured foliation} $F$ on $(F_g^s,\Delta)$
 is a foliation  which is transverse to $\Delta$ and which satisfies the following
four properties :

\medskip 

\noindent (1.4.1)\hskip .1in Each puncture of $F_g^s$ has  a 
neighborhood which is topologically 
an annulus and on which  
$F$ induces a foliation by    circles  
homotopic to the puncture.

\medskip 

\noindent (1.4.2)\hskip .1in 
The
foliation induced by $F$ on each  face of $\Delta$ has exactly one singular point, and the
  local model of that singular point  is  
  a 
    tripod, as in Figure 2  above.

\medskip 

\noindent (1.4.3)\hskip .1in The foliation induced by $F$ on each face of
$\Delta$ is 
equipped with a transverse measure, and  the total measure
 of
each transverse
arc  having 
 one endpoint at
 a vertex is infinite.

\medskip 

\noindent (1.4.4)\hskip .1in The transverse measures   induced on each edge
  of $\Delta$ from   its two sides differ by a homothety.

\medskip 

\noindent {\bf Equivalence relation.---} We shall say that two broken measured
foliations on  $(F_g^s,\Delta)$ are
equivalent if they differ by a transverse measure-preserving isotopy which preserves
 the
set of edges of $\Delta$.

\medskip

We denote by 
${\cal BM}(\Delta)$
 the  set of equivalence classes of 
 broken measured foliations relative to the  
   ideal triangulation $\Delta$.  Let ${\cal MF}$ denote the space of   equivalence classes
of measured foliations so that the each leaf has compact closure, or
equivalently, there are no leaves
running to the punctures. 
    We have a natural inclusion ${\cal MF}\subset {\cal BM}(\Delta)$. Measured foliations correspond to broken measured foliations for
    which the homothety factors in (1.4.4) are all equal to unity.
    Indeed, every measured foliation which 
    has no leaves running between punctures can be made transverse to   $\Delta$ by an isotopy, so that the subspace 
    ${\cal MF}(\Delta)$ of broken measured foliations with all the homothety factors equal to unity is canonically identified with ${\cal
MF}$.

    \medskip
    
There is a  natural action of the multiplicative group $\real_+=\{ t\in\real : t>0\}$ on 
${\cal BM}(\Delta)$,
obtained  by multiplying the transverse measure of each
broken measured foliation by a constant factor.
We denote the quotient space by ${\cal PBM}(\Delta)$.

\medskip

\noindent {\bf Definition 1.5 (Decorated broken measured foliation).---} A
{\it decorated broken measured foliation} on $F_g^s$ is a broken measured foliation
together with the choice, for each puncture of $F_g^s$,  of a closed leaf
homotopic to this puncture, in one of the annuli that are referred to  in
Property 1.4.1  of Definition 1.4.

 \medskip

 \noindent {\bf Example 1.6} Again, 
 a basic example of a decorated broken measured foliation is
 the horocyclic foliation associated to a decorated broken hyperbolic 
 structure. Here, the transverse measure for the induced foliation on each component
 of $F_g^s-\Delta$ is defined so as to coincide on each edge of that component
 with the Lebesgue measure induced from its structure as a hyperbolic ideal
 triangle.

  \medskip

  The space of decorated broken measured foliations, up to measure-preserving
  isotopy which fixes  the edges of $\Delta$, will be denoted by
  $ \widetilde{{\cal BM}}(\Delta)$,
  and the associated  projective space by $ \widetilde{{\cal PBM}}(\Delta)$.
Define the space $\widetilde{\cal MF}$ to be the collection of all measured foliations in ${\cal MF}$
together with a choice of embedded horocyclic leaf about each puncture and its corresponding projectivized
version $\widetilde{{\cal PMF}}$. 
   In the same way as for the non-decorated versions, we have natural inclusions
  $\widetilde{{\cal MF}} \subset \widetilde{{\cal BM}}(\Delta)$ and 
  $\widetilde{{\cal PMF}} \subset  \widetilde{{\cal PBM}}(\Delta)$.

\vskip .3in

\noindent  $\underline{ \bf
\hbox{2.--- Topology}}$ 
\vskip .2in

   A {\it triangle-edge
  pair} of $(F_g^s,\Delta)$ is a pair $(t,e)$ where $t$ is a triangle of
  $\Delta$ and $e$ an edge   of $t$.

   \medskip

We let $P$ be the set of triangle-edge pairs of $(F_g^s,\Delta)$.

  \medskip
  
 We define now the
 shift parameters   associated to the
   edge-pairs, which are  
   useful   parameters for the spaces ${\cal BH}(\Delta)$ and
   ${\cal BM}(\Delta)$.   The shift parameters for both spaces 
    are defined in similar manners, and we shall 
    use them to define the topology of
    these spaces.
    
    \medskip

    We start with the
   space ${\cal BH}(\Delta)$.

\medskip

   \noindent {\bf Shift parameters for 
    broken hyperbolic structures}.
    
     We first recall 
   that
   if $t$ is a hyperbolic ideal triangle, then each edge of $t$
   is equipped with a
   distinguished point, which is the foot of the perpendicular on that edge
   issuing from the center of   the triangle. Given 
    a broken
   hyperbolic structure $H$ on $(F_g^s,\Delta)$ and given an edge $e$ of $\Delta$, 
   there are
   two distinguished  points on   $e$,  
   corresponding to the
   inclusion of $e$ as an edge of  two triangles in $\Delta$,  each
    of these 
   triangles being 
   equipped with the structure of a hyperbolic ideal triangle inherited from
   $H$.

   Let $(t,e)$  be a triangle-edge pair 
   of $(F_g^s,\Delta)$.
     We   measure the distance between the  two distinguished 
    points on $e$ using
     the metric of  the triangle $t$. This distance is equipped with a
   {\it sign}, which is defined using  the convention given in Figure 3. 
   We note that this
   convention depends only on the orientation of the surface $F_g^s$, and not on
   any choice of 
   an orientation for the edge $e$. We denote this signed distance by $s_H(t,e)$, and
   we call it the {\it shift parameter} induced on  the triangle-edge pair
   $(t,e)$ by the structure $H$. We note that if $(t',e)$ is the other triangle-edge pair
   containing the edge $e$, then the shift parameter  $s_H(t',e) $ has the same
   sign as $s_H(t,e)$, and that  the values of the two parameters differ 
   by  a  multiplicative   factor
   equal to the   homothety factor in (1.1.1).

\medskip


\eject

~~\vskip 2.6in

\hskip .1in\epsffile{WP3ai.epsf}

\vskip .1in

\centerline{{\bf Figure 3} In case (a), the sign is positive, and in case
(b), the sign is negative.}

\medskip
   \medskip
   
   \noindent {\bf Shift parameters for hyperbolic structures}. 
   
   In the case where the broken
   hyperbolic structure $H$ is a hyperbolic structure, then the 
   shift parameters $s_H(t,e)$ and 
   $s_H(t',e)$ are equal, and we shall denote this common value by $s_H(e)$.

   \medskip
   
   \noindent {\bf Topology for broken hyperbolic structures}.
   
   It is clear that the  collection of shift parameters  associated to all
   the triangle-edge
   pairs of   $(F_g^s,\Delta)$ determines completely the broken hyperbolic
   structure. Thus, we have  an injective map from 
   ${\cal BH}(\Delta)$ into the space $\real^P$ of real valued functions on the
   set $P$  of triangle-edge
   pairs. Using this injection, we define  
     a topology for the
   space ${\cal BH}(\Delta)$ by taking the product topology on  
   $\real^P$ and the induced topology on the space ${\cal BH}(\Delta)$.
   Thus,  two broken hyperbolic structures
   are close if and only if the    shift parameters
   that they induce  on 
   triangle-edge pairs are close.

   \medskip

   We  now define  similar parameters for the space 
   ${\cal BM}(\Delta)$.

   \medskip

   \noindent {\bf Shift parameters for broken measured foliations}.
   
    Let $F$ be
   a broken measured foliation on $(F_g^s,\Delta)$
     and let $ e$  be again an edge of $\Delta$.
   On each side of $e$, there is a triangle of $\Delta$ equipped with a
   measured foliation induced from   $F$, there is a unique  singular
   point in the interior of this triangle, and there is a well-defined singular leaf issuing from that
   singular point and hitting the edge $e$. 
   The  two hitting points associated to the two sides of $e$ give two distinguished
   points on $e$. 
   If $(t,e)$ is now a triangle-edge pair in $(F_g^s,\Delta)$, then
   the segment in $e$ which joins these points has a well-defined measure,
 with respect to the transverse measure of the foliation induced on   the
   triangle $t$.  Again, we assign a
   {\it sign} to this measure, 
    using the convention given in Figure 4, and we note that the definition of this sign 
    depends  on the orientation of the surface $F_g^s$ and not on
a choice of an  orientation on $e$. We denote this signed measure by $s_F(t,e)$, and
   we call it the {\it shift parameter} induced on    the triangle-edge pair
   $(t,e)$ by the structure $F$. We note that if $(t',e)$ is the other triangle-edge pair 
   containing   $e$, then the quantity  $s_F(t',e) $  differs from
   $s_F(t,e)$ by a multiplicative factor which is equal
    to the   homothety factor in (1.4.4).

\medskip


~~\vskip 2.8in

\hskip .1in\epsffile{WP4ai.epsf}

\vskip .1in

\noindent {{\bf Figure 4}~In case (a), the sign is positive, and in case (b), the sign
is negative.}

\medskip

   \medskip
   
   \noindent {\bf Shift parameters for measured foliations}. 
   
   In the case where the broken
   measured foliation $F$ is a measured foliation, then  
   we have $s_F(t,e) =s _H(t',e)$,   and we  use the notation $s_F(e)$ 
   to denote this common value, which we call  the shift parameter
    associated to the edge $e$ by the measured foliation $F$.
   
   \medskip

   \noindent {\bf Topology for broken measured foliations}.
   
   The  collection of shifts associated to the various  triangle-edge
   pairs of $(F_g^s,\Delta)$ determines completely the broken measured
   foliation $F$. Thus, we have  an injective map from 
   ${\cal BM}(\Delta)$  to the space $\real^P$ of real valued functions on the
   set $P$ of triangle-edge pairs. This provides the  
   space ${\cal BM}(\Delta)$ with a topology
   for which two broken measured foliations
   are close if and only if 
    the  associated  set of shift parameters on all the  
   triangle-edge pairs are close.

   \medskip

  There is a natural map  ${\cal BH}(\Delta)\to {\cal BM}(\Delta)$
  which is defined 
  by associating to each broken hyperbolic structure its 
   broken
  measured foliation defined  in example 1.6 above and noting that, in that
  construction, the equivalence class
  of the broken measured foliation depends only on the equivalence class on the
  broken hyperbolic structure.

\medskip

Conversely, given a broken  measured foliation
  on $(F_g^s,\Delta)$, we can assign to it a broken hyperbolic structure by
  considering, for each 
    face of $
  \Delta$   with its induced foliation, a structure of  hyperbolic
  ideal triangle for which this measured foliation is the horocyclic foliation.
  Properties (1.4.2) and (1.4.3) insure  that this hyperbolic structure exists. 
  Now the
    structure of  hyperbolic ideal triangle on
   each 
   face of $
  \Delta$ makes the surface $F_g^s$ naturally equipped with a broken hyperbolic
  structures with the homothety factor  in (1.1.1)
  along each side  of $\Delta$ equal to the homothety
  factor  in (1.4.4). It
is clear that the equivalence class of this 
  broken hyperbolic structure depends only on the
  equivalence class of the broken measured foliation with which we
started.

  \medskip

  Thus, we have a map ${\cal BH}(\Delta)\to {\cal BM}(\Delta)$ which is one-to-one.
  We have the following proposition,  whose proof is clear from the description of the topologies on the spaces 
  ${\cal BH}(\Delta)$ and ${\cal BM}(\Delta)$.
  
  \medskip
  
  \noindent {\bf Proposition 2.1.---}{\it The map 
  ${\cal BH}(\Delta)\to {\cal BM}(\Delta)$  which assigns to the
  equivalence class of each broken hyperbolic structure  the equivalence class of the
  associated broken measured foliation is a homeomorphism.}
   \cqfd

  \medskip

  \noindent {\bf Topologies for 
  the decorated spaces.} 
  
  To define the topologies of the decorated spaces,
  we use the shift parameters that were defined  for the broken 
hyperbolic
  structures (respectively the broken measured foliations), together with extra
  parameters which give, for 
  each puncture of $F_g^s$, the distance between the distinguished closed curve
  corresponding to the decoration
  and a fixed point on an edge abutting on that puncture. In fact, we can 
  choose once and for all, for each puncture of $F_g^s$, one of the edges
   abutting on that puncture, and
  we fix a point on that edge. We then associate to 
  each decoration (that is, to  
  each distinguished closed 
  curve 
  around a puncture), the signed distances between the intersection points of
  this curve with the chosen edge  of $\Delta$  
  abutting on that puncture 
  and the    point that we fixed 
  on  that edge. To measure this distance,  we  choose one side for each edge,
   and we use the hyperbolic metric on the chosen 
  side. 
  The extra parameters for the decorated  broken measured foliations are defined in
  the same manner.
  
  \medskip

  It is clear now that  the space 
   $\widetilde{{\cal BH}}(\Delta)$ (respectively
    $\widetilde{{\cal BM}}(\Delta)$) 
  is an $\real^s$-bundle over the space ${\cal BH}(\Delta)$
  (respectively ${\cal BM}(\Delta)$), and we furthermore have the following

  \medskip

  \noindent {\bf Proposition 2.2.---}{\it The map 
  $f_{\Delta}: \widetilde{{\cal BH}}(\Delta)\to \widetilde{{\cal BM}}(\Delta)$ 
   which assigns to each
  equivalence class of decorated broken hyperbolic structure  the equivalence
   class of the
  associated decorated broken measured foliation is a homeomorphism.}\cqfd
  \medskip

\vskip .3in

\noindent  $\underline{ \bf
\hbox{3.--- The two-form on the space of decorated broken hyperbolic structures}}$ 
\vskip .2in

\noindent 
{\bf Definition 3.1 ($\lambda-$lengths).---}
Let $H$ be   an element  of $\widetilde{{\cal BH}}(\Delta)$
and let $t$ be a triangle  of
$\Delta$. Then $H$ equips  $t$ with the structure of a 
 hyperbolic  ideal triangle together  with three distinguished horocyclic
arcs centered at the vertices and joining pairwise the three edges of this
triangle, as in Figure 5. We note that the three  horocycles are
two-by-two disjoint since  by assumption they are leaves of the
horocyclic partial foliation of that triangle; this is in contrast to the situation in the
paper  [PP],  where the horocycles were allowed to intersect. For each edge $e$ of $t$, let
$\delta_H(t,e)$ be the
  distance between the two  horocycles which  intersect  $e$. (Thus, this distance is nonnegative, unlike the situation in [PP].)
 We define  
$\lambda_H(t,e)$,
 the {\it $\lambda-$length  of 
the triangle-edge pair $(t,e)$} with respect to  the broken
hyperbolic structure $H$, by the formula
$$\lambda_H(t,e)=\sqrt{2\exp(\delta_H(t,e))}.$$

\medskip


~~\vskip 2.2in

\medskip \hskip .9in\epsffile{WP5ai.epsf}

\vskip .1in

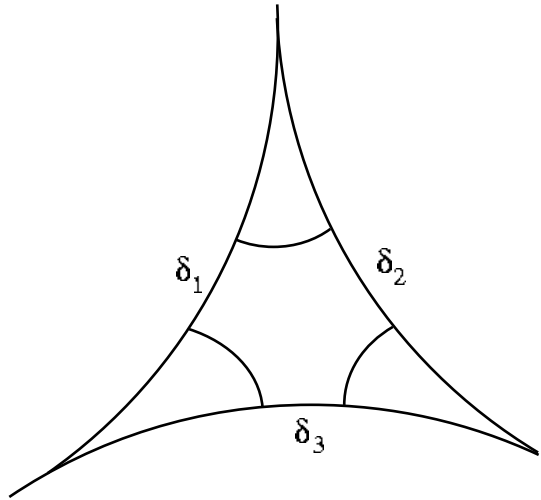
\centerline{{\bf Figure 5}~The the three horocycles induced by the
decoration.}

\medskip

\medskip

\noindent {\bf Proposition 3.2} {\it A decorated broken hyperbolic structure is uniquely
determined by its collection of
$\lambda-$lengths of triangle-edge pairs.}
  
  \medskip
  
  \dem
   This can be seen
at the level of a face of $\Delta$: Let $t$ be a
hyperbolic ideal triangle  with vertices  $v_1$, $v_2$, $v_3$, equipped with its 
 three horocyclic arcs joining the edges pairwise
 as in Figure 5, with relative distances  $\delta_1$, $\delta_2$ and $\delta_3$.
The proof of Proposition 3.2 is a consequence of the following

\medskip

\noindent {\bf Lemma 3.3.---} 
{\it The three numbers $\delta_1$, $\delta_2$ and $\delta_3$ 
  determine completely  the positions of the three horocyclic arcs in the
triangle
$t$.}

\medskip

\noindent \dem 
We show that it is not possible to  keep the same   distances between 
the horocyclic arcs if we move any of these horocycles from its original
position.
Suppose for that  that we move $h_1$ towards $v_1$. If we want to keep the 
same distance $\delta_2$ between $h_1$ and $h_3$,  then we have to move 
$h_3$ to a  new position  that
has to be farther away from the vertex $v_3$
than the old position. Now, in order to keep the distance $\delta_1$ between $h_3$ and $h_2$ unchanged,
 we have to move $h_2$ to 	a new position which is
nearer than the old one to the vertex $v_2$. But then, the new distance between $h_2$ and $h_1$ will be greater than $\delta_3$.
Thus, we cannot move $h_1$ towards $v_1$.
The same kind of argument shows that we cannot move $h_1$ away from $v_1$ while keeping unchanged  the three relative distances 
$\delta_1$, $\delta_2$, $\delta_3$ between  the horocyclic arcs.
This proves Lemma 3.3 and Proposition 3.2.\cqfd

\medskip

Thus, the map from $\widetilde{{\cal BH}}(\Delta)$ to the set  of functions on the collection
$P$ of triangle-edge pairs of $\Delta$, which assigns  to each
triangle-edge pair its
$\lambda$-length,  is injective. It is clear that this map is continuous
and proper   (the topology on
$\widetilde{{\cal BH}}(\Delta)$ being given by the shift parameters). Thus, we have the following

\medskip

\noindent {\bf Proposition 3.4.---} 
\it The map 
$\widetilde{{\cal BH}}(\Delta)\to \{ t\in\real _+:t\geq\sqrt{2}\}^P$
induced by $\lambda$-lengths is a homeomorphism onto its image.\rm

\medskip

\noindent 
{\bf Definition 3.5 (The two-form ${\widetilde \Omega}$ on $\widetilde{{\cal BH}}(\Delta)$).---}
We  define the  form ${\widetilde \Omega}$ on  
$\widetilde{{\cal BH}}(\Delta)$  in terms of the $\lambda-$length  coordinates by the formula 
$${\widetilde \Omega} = -2\sum d\log a\wedge d \log b + 
d\log b\wedge d\log c + 
d\log c\wedge d\log a,$$
where the sum is taken over   the set of   triangles  of
  $\Delta$, and where 
   for each such  triangle,  $a$, $b$ and $c$  are 
  the $\lambda-$lengths of the three associated 
  triangle-edge pairs taken   in a counterclockwise
  cyclic order. 
  
From the way it is defined, it is clear that $\widetilde\Omega$ 
 on
  $\widetilde{{\cal BH}}(\Delta)$ pulls back to the (restriction to $\widetilde{\cal T}'\subseteq\widetilde{\cal T}$ of) 
form 
  ${\widetilde  \omega}$ 
   which was defined
    in [Pe2] on the decorated Teichm\"uller space $\widetilde{{\cal
T}}\subset\widetilde{{\cal BH}}(\Delta)$ of the surface
    and which we used
      in the paper [PP]. The form  ${\widetilde  \omega}$, in turn, is the pull-back of the
Weil-Petersson K\"ahler two-form on the
      Teichm\"uller space ${\cal T}$ of the surface; in this sense $\tilde\Omega$ ``extends''
the Weil-Petersson two-form.

\vskip .3in

\noindent  $\underline{ \bf
\hbox{4.--- The two-form on the space of decorated broken measured foliations}}$ 
\vskip .2in

\noindent {\bf Definition 4.1 (Compactly supported broken measured foliation).---}
Consider an element of ${\cal BM}(\Delta)$  represented by a foliation $F$ on $F_g^s$.   
We  associate to $F$ a partial foliation $F_0$ on $F_g^s$ by deleting in the neighborhood of 
 each puncture of $F_g^s$  the largest annulus foliated by closed leaves of $F$ which is
described in
Property (1.4.1) of Definition 1.4. We call this foliation $F_0$ the {\it compactly supported broken measured foliation on $F_g^s$ associated to $F$}.

\medskip

We regard two  compactly supported broken measured foliations as equivalent if 
they are associated to two broken measured foliations representing the same element  of
${\cal BM}(\Delta)$, or, equivalently, if one of these 
compactly supported broken measured 
foliations can be obtained from the other one by an isotopy of the surface which preserves setwise each edge
of $\Delta$ and the transverse measure in each face of $\Delta$.

\medskip

We denote by
${\cal BM}_0(\Delta)$ the set of equivalence classes of compactly supported broken measured 
foliations on $F_g^s$.

\medskip

We let ${\cal MF}_0$ be the set of equivalence classes of compactly supported measured foliations on $F_g^s$. 
Again, we note that any element of ${\cal MF}_0$ can be represented by   a compactly supported 
 broken 
measured foliation on $(F_g^s,\Delta)$  having all of its homothety factors  equal to unity. Thus, 
there is a natural inclusion
${\cal MF}_0\subset {\cal BM}_0(\Delta)$, for each  ideal triangulation $\Delta$. 
Furthermore, it is clear  that the construction of the compactly supported
broken measured foliation in Definition 4.1 gives
 a natural ``identification'' ${\cal BM}(\Delta)
\buildrel\approx\over\to
{\cal BM}_0(\Delta)$ (i.e., a homotopy inverse to the inclusion ${\cal BM}_0(\Delta)
\subset
{\cal BM}(\Delta)$) 
which restricts to an identification between the two subspaces ${\cal MF}\subset {\cal BM}(\Delta )$ and ${\cal
MF}_0\subset {\cal BM}_0(\Delta)$.

\medskip

\noindent {\bf Definition 4.2 (Compactly supported broken measured foliation with collars)} 
A {\it compactly supported broken measured foliation with collars} is a partial broken measured foliation on $F_g^s$ 
which is obtained from a broken measured foliation
by removing a (not necessarily maximal)   annulus neighborhood of each puncture
of $F_g^s$, which is foliated by circles. As in the case of compactly supported broken measured
foliations, there is a natural equivalence relation on the space of 
compactly supported broken measured foliations with collars,  and we let 
 $\widetilde{{\cal BM}}_0(\Delta)$ be the space of equivalence classes.

 \medskip
 
 We let $\widetilde{{\cal PBM}}_0(\Delta)$ be the space of equivalence classes of compactly supported broken measured foliations with collars, 
 up to the natural  
 action of $\real_+$ on the transverse measures.
 
 \medskip

 To each element of  $\widetilde{{\cal BM}}(\Delta)$, we   associate  an element of 
 $\widetilde{{\cal BM}}_0(\Delta)$, by removing,  from the support
 of  a broken measured foliation $F$ 
 representing the element of
 $\widetilde{{\cal BM}}(\Delta)$,
 the foliated annulus around each puncture of
$F_g^s$ which is bounded by the closed leaf corresponding to the decoration.
It is clear that this gives a natural identification
$ \widetilde{{\cal BM}}(\Delta) 
\buildrel\approx\over\to  \widetilde{{\cal BM}}_0(\Delta)$.

  \medskip
  
 We let  ${\widetilde{\cal MF}_0}'$
be the set of equivalence classes of compactly supported measured foliations with collars on $F_g^s$. 
Note that there is a natural one-to-one correspondence between the space  ${\widetilde{\cal
MF}_0}'$ and the space 
 $\widetilde{\cal MF} $ which is defined at the end of \S1 above.
 Again, we note that there is a natural inclusion
 $\widetilde{{\cal MF}_0}'\subset \widetilde{{\cal BM}}_0(\Delta)$,  
 for every ideal triangulation $\Delta$.
Recall from [PP] the space $\widetilde{\cal MF}_0$ of decorated measured foliations of compact support; in that work, we considered formal
collar widths, which may be positive or negative, where the space of collars of non-negative widths is identified with
${\widetilde{\cal MF}_0}'$. 

\medskip

The identification
$ \widetilde{{\cal BM}}(\Delta) 
\buildrel\approx\over\to  \widetilde{{\cal BM}}_0(\Delta)$  restricts to an identification between the subspace 
$\widetilde{{\cal MF}}\subset \widetilde{{\cal BM}}(\Delta) $ of decorated measured foliations 
and the subspace $\widetilde{{\cal MF}_0 }'\subset 
\widetilde{{\cal BM}}_0(\Delta)$  of compactly supported   measured foliations
with collars  on  $F_g^s$.

\medskip

The next proposition follows directly  from the definitions.

\medskip

\noindent {\bf Proposition 4.3} {\it  There is a natural fibre bundle 
$\widetilde{{\cal BM}}_0(\Delta)\to {\cal BM}_0(\Delta)$  whose fibre above a point is $(\real_+\cup\{ 0\} )^s$ (the set of collar weights),
which restricts to a fibre bundle
 $\widetilde{{\cal MF}_0}'\to {\cal MF}_0$ over  the subspace
 $\ {\cal MF}_0$ of ${\cal BM}_0(\Delta)$.}

\medskip

\noindent 
{\bf  The dual graph $G$.}
The dual graph  to the triangulation $\Delta$  is  
a graph $G$ which is embedded in $F_g^s$ and which is 
defined as follows. In each face of $\Delta$, there is a vertex of
$G$, and   two such vertices are joined    by a segment in the graph
(and this segment will be the union of two edges of the graph)   whenever the
corresponding faces  share a common  edge. The segment intersects 
the corresponding edge  of $\Delta$ in a unique 
point,  and 
we take this point to be  also a vertex
of $G$. 
Thus, the graph $G$ has two kinds of vertices: trivalent vertices (which are
in one-to-one correspondence with the faces of the triangulation $\Delta$), and bivalent
ones (in one-to-one correspondence with the edges of $\Delta$). 
An {\it  edge} of $G$ is then a connected component of the space $G-\{
\hbox{vertices}\}$. 

\medskip

We note 
by the way that the graph 
$G$ together with its embedding in $F_g^s$ is an example of what is
 usually called a {\it fatgraph} (that is, a graph equipped with
a cyclic ordering at each of its vertices).

\medskip

We shall define  a two-form on $\widetilde{{\cal BM}}(\Delta)$
which is an extension of the form $\iota$  which we defined in ([PP], \S3) on  
the space $\widetilde{{\cal MF}}$ of decorated measured foliation in $F_g^s$ and which in turn
was an extension of Thurston's form on the space of measured foliations on $F_g^s$.
For the purpose of making this definition, we need to consider, as we did   in ([PP], \S 2), 
  the null-gon track or ``freeway'' which is dual to $\Delta$. 

\medskip

\noindent  {\bf  The freeway dual to $\Delta$.}  The freeway $\tau$ dual to $\Delta$ is a maximal train track in
$F_g^s$ in the usual sense except that the connected components of $F_g^s - \tau$ are allowed to be
 once-punctured monogons. In fact, there is one such component around  
each puncture of $F_g^s$. The
intersection of $\tau$  with each
face $t$ of $\Delta$ is as in Figure 6. Thus,  there are three
trivalent  vertices of $\tau$ in the interior of
 $t$. Furthermore, there is   a 
bivalent vertex of
$\tau$ on each edge of $\Delta$, in the same way 
 as for the  graph $G$ dual to $\Delta$. In fact, $\tau$ is obtained from $G$ by blowing up each trivalent vertex of 
$G$ into a triangle as illustrated in Figure 6.

\medskip


~~\vskip 2.2in

\hskip .9in\epsffile{WP6ai.epsf}

\vskip .1in

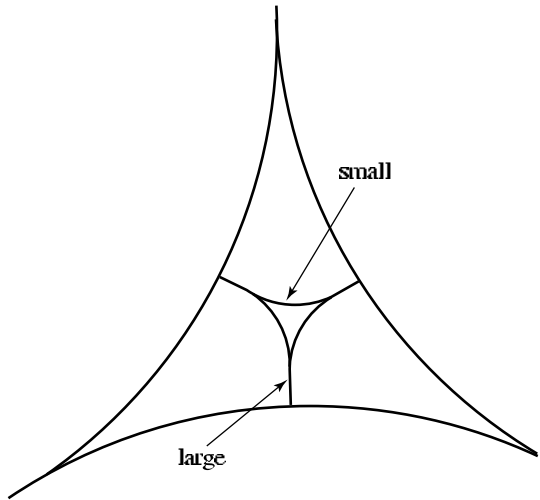
\centerline{{\bf Figure 6} The intersection of $\tau$ with a face of
$\Delta$.}

\medskip

It follows from these definitions
that  each edge of $\tau$ is contained in a face of $\Delta$, 
and that there are two types of edges in each such   face: the edges  which have one vertex on 
the boundary of that face, which we call the {\it large} edges, and those 
whose two vertices are in the interior of
that face, which we call the {\it small} edges (see Figure 7).

\medskip

We define now
a system of 
weights on the edges of a
$\tau$, associated to a  
 compactly supported foliation $F_0$.  
 We deal first  with 
the large edges. Let $e$ be   a large edge of $\tau$.
 Then, $e$
 is contained in a triangle $t$ of $\Delta$, and $e$ has 
    one of its endpoints on an edge $\ell$ of $\Delta$. The weight $w(e)$ induced by $F_0$ on $e$ 
is defined as  the total  measure of the edge $\ell$ with respect to the
transverse measure of the restriction of the foliation $F_0$ to the triangle
 $t$. Consider now a face $t$ of
 $\Delta$, and let $w(a)$, $w(b)$ and $w(c)$ be the weights
which are defined in this way on the three large edges $a$, $b$ and $c$
of $\tau$ that are contained 
 in  $t$. There are   three   weights, 
 $w(\alpha)$, $w(\beta)$ and $w(\gamma)$ 
induced on the three small  edges 
$\alpha$, $\beta$ and $\gamma$ contained in 
    $t$. They  are defined by the three   equations 
    
    $$ (4.1.1)\hskip .4in \cases{w(\alpha) =(1/2)~\bigl [(w(b)+w(c)-w(a)\bigr ],&\cr\cr
 w(\beta)=(1/2)~\bigl [w(c)+w(a)-w(b)\bigr ],&\cr\cr
 w(\gamma)=(1/2)~\bigl [w(a)+w(b)-w(c)\bigr ].\cr}$$
 
  \medskip

 We note that the triple $w(a)$, $w(b)$, $w(c)$ satisfies the three large triangle inequalities, as a consequence of  the 
 condition of conservation of mass
 at the trivalent vertices of $\tau$, and therefore the three weights $w(a)$, $w(b)$ and $w(c)$  are $\geq 0$.

\medskip


~~\vskip 2in

\hskip .9in\epsffile{WP7ai.epsf}

\vskip .1in

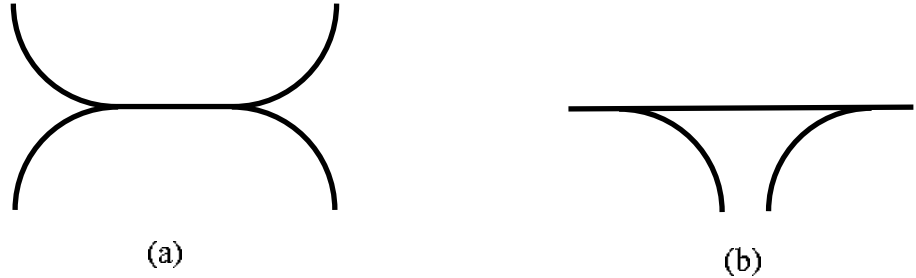
\centerline{{\bf Figure 7} Case (a) is a large edge and case (b) a small
edge.}

  \medskip
 
Any system of nonnegative weights on the edges of $\tau$ gives a well-defined element of the space 
 ${\cal BM}_0(\Delta)$ of equivalence classes of compactly supported broken measured foliations 
 with collars on $F_g^s$, provided these weights satisfy the following   condition:

 \medskip
 
 \noindent (4.1.2) \hskip .1in
  At each trivalent vertex of $\tau$, there is a condition  of conservation
 of mass, as in usual train track theory: if at this trivalent vertex, 
  $\alpha$ and $\beta$ are the two edges abutting from one side, and $c$ the edge abutting 
  from the other side, and if  $w(\alpha)$, $w(\beta)$ and $w(c)$ are the
 weights on these three edges, then we have $w(\alpha)+w(\beta)=w(c)$.

 \medskip

\noindent {\bf Definition 4.4 (Broken measure on the dual freeway).---} We shall say that a system of nonnegative weights on the edges 
 of $\tau$ satisfying   condition (4.1.2) on the conservation of mass a {\it broken measure} on $\tau$.
  
  \medskip

We note again  that a broken measure  
  on $\tau$ can be  obtained by taking an arbitrary system  of nonnegative
  weights on the large edges  of $\tau$  satisfying the three large triangle inequalities and condition (4.1.2),
since such a system of weights induces a unique system of weights (on all the edges of $\tau$).

 \medskip

 To each element of $B(\tau)$ which is not identically zero, we associate a partial 
 foliation by the usual construction which replaces each edge of $\tau$ having  nonzero weight by a rectangle foliated by
 leaves parallel to its 
  ``horizontal'' sides, 
these sides being  nearly parallel to that edge. Each foliated rectangle is equipped with
 a transverse measure whose total mass is equal to the weight of the corresponding edge. 
 The various rectangles are then glued along their ``vertical'' sides by measure-preserving maps at the trivalent vertices of $\tau$ and by affine maps at the
 bivalent ones, and we obtain a partial foliation on $F_g^s$.
 It is clear that this partial foliation is an element of
  $\widetilde{{\cal BM}}_0(\Delta)$. Thus, we have a natural map 
  which assigns to each  element of 
  $B(\tau)$  which is not
  non-identically zero  an element of $\widetilde{{\cal BM}}_0(\Delta)$. We assign  to the zero element of
   $B(\tau)$  the empty foliation, and we thus obtain a natural map
   $  B(\tau)\to \widetilde{{\cal BM}}_0(\Delta)$ (where we have included the empty foliation as an element of
   $\widetilde{{\cal BM}}_0(\Delta)$).
   We recall now that there is a natural identification between the spaces 
    $\widetilde{{\cal BM}}_0(\Delta)$
   of compactly supported broken measured foliations with collars	and the space
    $\widetilde{{\cal BM}}(\Delta)$
   of decorated broken measured foliations. Composing the two maps, we obtain a map 
   $h_{\Delta}: B(\tau)\to \widetilde{{\cal BM}}(\Delta)$ . 
We have the  following

\medskip

\noindent {\bf Proposition 4.5.---} {\it The natural map $h_{\Delta}: B(\tau)\to \widetilde{{\cal BM}}(\Delta)$ is a homeomorphism onto.}

\medskip

\dem The only non-trivial fact is the injectivity of this map, and this is can be seen in the same way as in the proof of Theorem 
3.1 of [Pa], where the same injectivity property is proved in the case of closed surfaces.

 \medskip

 Proposition 4.5 provides a set of global coordinates for the space
 $\widetilde{{\cal BM}}(\Delta)$.
 \medskip

\noindent 
{\bf Definition 4.6 (The two-form ${\widetilde \iota}$ on $\widetilde{{\cal BM}}_0(\Delta))$.---}
Using the coordinates  provided by proposition 4.5 on $\widetilde{{\cal BM}}(\Delta)$,  we define the two form 
  ${\widetilde \iota}$ on this space by the formula
  $${\widetilde \iota} = -(1/2)\sum d w(\alpha) \wedge d w(\beta) + 
d w(\beta) \wedge  d w(\gamma) + 
d w(\gamma) \wedge d w(\alpha),$$ 
where the sum is taken over   the set of   triangles   of
  $\Delta$, and where 
   for each such triangle,  $\alpha$, $\beta$ and $\gamma$  are 
  the small edges of $\tau$  
  which are contained in this triangle, in counterclockwise order with respect to the orientation of the surface,
  and  where $w(\alpha)$, $w(\beta)$ and $w(\gamma)$ are the weights on these edges which are induced by the element of
  $\widetilde{{\cal BM}}_0(\Delta)$.
 By the natural homeomorphism
 $\widetilde{{\cal BM}}_0(\Delta)\to \widetilde{{\cal BM}}(\Delta)$, we can consider the form 
 ${\widetilde \iota}$ as being defined on the space
 $\widetilde{{\cal BM}}(\Delta)$
  
  \medskip

 It is clear from the definitions that this form is an extension of the form  $\iota$ defined in ([PP], \S3) on  
  the space $\widetilde{{\cal MF}_0}$   of decorated 
  compactly supported   measured foliations on $F_g^s$, which
  in turn is an extension of Thurston's symplectic form on the space ${\cal MF}_0$ of compactly supported measured foliations on $F_g^s$.

\vskip .3in

\noindent  $\underline{ \bf
\hbox{5.--- The relation between the various two-forms}}$ 
\vskip .2in

  By 
   Proposition 2.3, we have a natural homeomorphism
  $ f_{\Delta}:\widetilde{{\cal BH}}(\Delta)\to \widetilde{{\cal BM}}(\Delta)$.

   \medskip

 \noindent {\bf   Proposition 5.1.---}
 {\it   The  homeomorphism
$f_{\Delta}:\widetilde{{\cal BH}}(\Delta)\to \widetilde{{\cal BM}}(\Delta)$ preserves the forms
${\widetilde \Omega}$ on $\widetilde{{\cal BH}}(\Delta)$
and
${\widetilde \iota}$ on $\widetilde{{\cal BM}}(\Delta)$.} 

\medskip

\dem Of course, we use the $\lambda-$length
coordinates on $\Delta$ for the space 
$\widetilde{{\cal BH}}(\Delta)$ and the broken measures on the dual freeway $\tau$ for the space
$\widetilde{{\cal BM}} (\Delta)$.
Since the two-forms are defined in both  cases by taking sums over triangles of $\Delta$, it suffices to consider the contributions of the forms to each
such  triangle.

Let  $H\in  \widetilde{{\cal BH}}(\Delta)$.
We   describe the transformation $f_{\Delta}$ on the coordinates associated to  $t$. Consider the dual graph $G$. 
Since each edge of $G$ is contained in a unique
triangle-edge pair of $(F_g^s,\Delta)$,
we  can consider   the 
  $\lambda-$lengths associated to $H$ as being  defined
on the edges of the graph $G$, and we denote $\lambda-$length of such an edge $e$ by $\lambda(e)$.
 To each such   edge $e$  is naturally associated a large edge $e'$ of the dual
 freeway
 track $\tau$. By examining the  way the homeomorphism of Proposition 2.3 is defined, we can see that the weight $w(e')$ on $e'$ 
 which is induced by the transformation $f_{\Delta}$ is given by
  $$w(e')=2\log \lambda(e)+\log (1/2).$$
 Taking differentials, we have $d w(e')=2d\log\lambda(e).$
 
 \medskip

  We now use equations (4.1.1) which relate the three small edges of the dual null-track $\tau$ which are in $\Delta$ 
  to the weights on the large edges, and  the formula for the two-form 
  ${\widetilde \iota}$ given  in Definition 4.6.
  This tells us that the contribution of
 ${\widetilde \iota}$   to the triangle $t$ of $\Delta$, 
 in which the three edges of the dual graph $G$  in counterclockwise order are denoted by $a$, $b$ and $c$, is equal to $-(1/2)$ times the following quantity:
 $$\bigl[  d\log \lambda(b) +  d \log \lambda(c) -  d \log \lambda(a)\bigr ]\wedge
\bigl[   d\log \lambda(c) +  d \log \lambda(a) -  d \log \lambda(b)\bigr ]$$
 $$+\bigl[  d\log \lambda(c) +  d \log \lambda(a) -  d \log \lambda(b)\bigr ]\wedge
 \bigl[  d\log \lambda(a) +  d \log \lambda(b) -  d \log \lambda(c)\bigr ]$$
 $$+ \bigl[ d\log \lambda(a) +  d \log \lambda(b) -  d \log \lambda(c)\bigr ]\wedge
 \bigl[  d\log \lambda(b) +  d \log \lambda(c) -  d \log \lambda(a)\bigr ].$$

 Simplifying this formula, we find that the 
  contribution of
 ${\widetilde \iota}$   to the triangle $t$ is equal to
 $$-2    d\log \lambda(b)\wedge d\log\lambda(c) + 2 d \log \lambda(c)
 \wedge d\log\lambda(a) - 2 \log \lambda(a)\wedge d\log\lambda(b).$$

Comparing with the formula 
for   the two-form
 ${\widetilde \Omega}$ given in Definition  3.5,
  we see that the last  quantity  is also equal to the 
  contribution  of 
 ${\widetilde \Omega}$   to the triangle $t$. This proves Proposition 5.1.\cqfd

\medskip

We consider now  the product ``Yamabe'' space 
${\widetilde{\cal Y}}={\widetilde{\cal BH}}(\Delta)\times ]0,\infty[$
in which  we view each element  as the equivalence class of a
decorated broken  structure obtained by gluing ideal triangles,   
whith the ideal triangles having
 constant curvature, but not necessarily equal to $-1$. 
More precisely, if $H$ is the equivalence class of a decorated broken
hyperbolic structure on $F_g^s$  and if $x\in ]0,\infty[$, then the
element $(H,x)$ of ${\widetilde{\cal Y}}$ denotes the equivalence class
(for the relation of isotopy   preserving the vertices, the edges  and the faces 
of $\Delta$) of a structure on $F_g^s$ obtained by replacing
each hyperbolic ideal triangle by an ideal triangle in the complete simply
connected Riemann surface  with Gaussian curvature equal to $-x^2$, 
and gluing the new triangles using the same $\lambda-$lengths as for the structure $H$.

 \medskip

We let $q:{\widetilde{\cal Y}}\to {\widetilde{\cal BH}}(\Delta)$ denote
the projection map onto the first factor.

\medskip

There is a natural topology on the union
${\overline{\cal Y}}={\widetilde{\cal Y}}\cup {\widetilde{\cal
BM}}_0(\Delta)$
which makes ${\overline{\cal Y}}$ homeomorphic to the space 
${\widetilde{\cal
BM}}_0(\Delta)\times [0,\infty[$  and 
which extends
the topology which we considered in \S5 of [PP] on the union
${\widetilde{\cal T}}\cup {\widetilde{\cal MF}}_0(\Delta)$.
This topology on  
${\overline{\cal Y}}$
is defined as in [PP], \S5, using the homeomorphism 
$f_{\Delta}:
{\widetilde{\cal BH}}(\Delta)\to {\widetilde{\cal BM}}(\Delta)$ of Proposition 2.3, instead of the map
${\cal F}_{\Delta}$ that we used in [PP]. The   topology 
on
${\overline{\cal Y}}$ is defined  from the usual Thurston techniques for the
compactification of Teichm\"uller space, and it has the property
  that the quotient of the space ${\overline{\cal Y}}$  by the
action of $\real_+$ (by homotheties on the second factor of the space   
${\widetilde{\cal Y}}={\widetilde{\cal BH}}(\Delta)\times ]0,\infty[$
and by multiplying the transverse measure by a constant factor  on the space
${\widetilde{\cal
BM}}_0(\Delta)$) gives a topology on the union 
${\widetilde{\cal BH}}(\Delta)\cup 
{\widetilde{\cal
PBM}}_0(\Delta)$
which is a natural extension  
of the topology on
${\cal BH}(\Delta)\cup 
{\cal
PBM}_0(\Delta)$
which was defined in [OP2]. 

\medskip

There is a criterion, which is equivalent to a criterion  stated as Theorem 5.5 in
 [OP1], for the convergence of a sequence of covering hyperbolic structures (see appendix B) to the projective
class of an affine foliation.   The criterion  in [OP1] uses the notion of covering hyperbolic structures instead of the
broken hyperbolic structures that we use here, and it uses affine foliations instead of the broken measured foliation that we use here.
The convergence criterion that we state below is tailored to our setting of broken hyperbolic structures and broken measured foliations, and
it is valid because the topology of the union of the two spaces, as it is defined in [OP2], \S7, like the topology of the corresponding
spaces in [OP1], as an endpoint compactification along each ray.
We shall not repeat the details of the proof, which are a straightforward  adaptation of those of [OP1].

\medskip

In the setting of this paper, the convergence criterion   
relies upon the map $f_{\Delta}:
{\widetilde{\cal BH}}(\Delta)\to {\widetilde{\cal BM}}(\Delta)$ of Proposition 2.3, which, by the canonical identifications between the
spaces ${\widetilde{\cal BM}}(\Delta)$ and ${\widetilde{\cal BM}}_0(\Delta)$, we can consider as a map
$$f_{\Delta}: {\widetilde{\cal BH}}(\Delta)\to {\widetilde{\cal BM}}_0(\Delta).$$

\medskip

\noindent {\bf Convergence criterion:} Let $(H_n,x_n)$ be a sequence of elements in  the space 
${\widetilde{\cal Y}}={\widetilde{\cal BH}}(\Delta)\times ]0,\infty[$ such that $H_n$ eventually leaves every compact set in
${\widetilde{\cal BH}}(\Delta)$, and such that $x_n\to 0$. Then, $(H_n,x_n)$ converges to an element $F\in {\widetilde{\cal BM}}_0(\Delta)$
in the topology of ${\overline{\cal Y}}$ if and only if the sequence $x_nf_{\Delta}(H_n)$ converges to $F$ in the topology of
 ${\widetilde{\cal BM}}_0(\Delta)$. 
 
 \medskip
 
 We define now  the map
 $$ {\widetilde h}: {\widetilde{\cal Y}}\to 
{\widetilde{\cal BM}}_0(\Delta)\times [0,\infty[$$ by 
the formula
$$(H,x)\mapsto (x f_{\Delta}(H) , x)$$
for every $(H,x)\in {\widetilde{\cal BH}}(\Delta)\times ]0,\infty[$.

\medskip

The image of the map ${\widetilde h}$ is the subspace
${\widetilde{\cal BM}}_0(\Delta)\times ]0,\infty[$. 

\medskip

By the convergence criterion above, this map extends to a homeomorphism
$$ {\overline h}: {\widetilde{\cal Y}}\to 
{\widetilde{\cal BM}}_0(\Delta)\times [0,\infty[$$
in which the image of the set $ {\overline {\cal Y}}- {\cal Y}$ of ideal points  is the subset 
${\widetilde{\cal BM}}_0(\Delta)\times \{0\}$.
We consider the two-form ${\widetilde \Omega}$ on
${\widetilde{\cal BH}}(\Delta)$ and its pull-back $q^*({\widetilde \Omega})$ on $ \widetilde{ {\cal Y}}$, which  we normalize to a form
$\Omega'=x^{-2} q^*({\widetilde \Omega})$ on ${\widetilde {\cal Y}}$.
Likewise, we let $\iota'$ be the pull-back of the two-form ${\widetilde \iota}$ on $\widetilde{{\cal BM}}_0(\Delta)$ by the projection of
$\widetilde{{\cal BM}}_0(\Delta))\times ]0,\infty[$ onto the first factor. Then, using  Proposition 5.1, $\Omega'$ extends continuously to a two-form
${\widetilde \Omega'}$ on
${\overline {\cal Y}}$
 whose restriction on
 $ {\overline {\cal Y}}- {\cal Y} \approx {\widetilde{\cal BM}}_0(\Delta)$ is the form ${\hat \iota}$ which extends Thurston's form on
 ${\widetilde{\cal BF}}_0(\Delta)$. 
 
 \medskip

 We summarize these results as
 
 \medskip

  \noindent {\bf Theorem 5.2.---}{\it The Weil-Petersson K\"ahler form on the  Teichm\"uller space ${\cal T}$ of $F_g^s$ extends to a
two-form 
   ${\widetilde \Omega}$ on the space $\widetilde{{\cal BH}}(\Delta)$ of decorated broken hyperbolic structures on $(F_g^s,\Delta)$,
   which induces a two-form $\Omega'$ on the space ${\widetilde{\cal Y}}={\widetilde{\cal BH}}(\Delta)\times ]0,\infty[$, which extends
continuously
    to a  two-form 
   ${\widetilde \Omega'}$ on
${\overline {\cal Y}}\cup {\widetilde{\cal BM}}_0(\Delta)$. The restriction of this form 
 ${\widetilde \Omega'}$  to the space ${\widetilde{\cal BM}}_0(\Delta)$ is an extension of Thurston's symplectic form on
the space
  ${\cal MF}_0$  of compactly supported
 measured foliations on the surface.}
 \cqfd

\medskip

It seems to us a basic question as to whether $\widetilde\Omega$ is non-degenerate (as are both the Weil-Petersson and Thurston forms).

\medskip 

We finally remark that one of the innovations in [PP] which has subsequently proved useful is the PL isomorphism of $\widetilde{\cal
MF}_0$ with the {\sl vector space} of measures on the freeway, which requires consideration of both positive and negative collars.  In the
current paper, we consider only positive collars in $\widetilde{\cal MF}_0$ and have the suspicion that there may be a nice synthesis of
broken structures with negative collars.

\vskip .3in

\noindent $\underline {\rm Appendix~A~-~Affine~foliations~and~broken~measured~foliations.}$

\vskip .2in

For completeness, we recall the definition of affine foliation from [OP1], which was adapted from [HO] to 
the case of punctured surfaces, and compare with broken measured foliations. 

\medskip

\noindent  {\bf Definition A.1.---(Affine foliation)} A
 foliation $F$ on   $F_g^s$
is   an {\it affine foliation} if the following four properties are satisfied :

\medskip

\noindent (A.1.1) \hskip .1in There is 
 a neighborhood  of each puncture of $F_g^s$ which 
is topologically an annulus on which $F$ induces a foliation by
 closed leaves which are homotopic to the puncture.

\medskip

\noindent (A.1.2) \hskip .1in  The foliation $F$ has no Reeb components.

\medskip

\noindent (A.1.3)  \hskip .1in The lift $\widetilde F$ of $F$ to
 the universal cover ${\widetilde F_g^s}$ is equipped
with a transverse measure $\mu$ such that for each covering translation
$\alpha$, there exists a positive real number
$\phi(\alpha)$ with the  property that for every arc $c$ in ${\widetilde F_g^s}$ 
which
is transverse to $\widetilde F$, we have $\mu(\alpha(c))=
\phi(\alpha)\mu(c)$.

\medskip

\noindent (A.1.4)  \hskip .1in If  $c$ is an arc 
which is transverse to  $F$ and which has one of its endpoints at a
puncture of  $F_g^s$, then the 
$\mu$-transverse measure of each lift of $c$ to ${\widetilde F_g^s}$  is 
infinite.

\medskip

If $\Gamma$ denotes  the group of covering translations of the cover
 ${\widetilde F_g^s}\to
F_g^s$, then  it is easy to see that the map $\alpha\mapsto \phi(\alpha)$ from
$\Gamma$ to $\real_+^*$ is a homomorphism. It is called the {\it holonomy
homomorphism} of the affine foliation. Using the canonical isomorphism
$\Gamma\simeq\pi_1(S)$, we can thus associate to each affine foliation on 
$F_g^s$   a
homomorphism $\pi_1(S)\to\real_+$, which is also called the holonomy
homomorphism.

\medskip

\noindent {\bf Equivalence relation.---} We say that  two affine foliations $F$
  and $F'$ on $F_g^s$ 
  are equivalent if $F$ can be obtained from $F'$ by an isotopy 
whose lift   to ${\widetilde F_g^s}$ preserves the
transverse measure of the lifted foliations ${\widetilde F}$ and ${\widetilde 
F'}$. 

\medskip

The set of equivalence classes of affine foliations is
denoted by
${\cal AF}$ and was first studied by Hatcher and Oertel in [HO].

\medskip 

Proposition 7.1 of [OP2] establishes a natural homeomorphism between ${\cal BM}(\Delta )$ and ${\cal AF}$, for any ideal
triangulation $\Delta$ of $F_g^s$. To define this homeomorphism, choose a face $F_1$ of $\widetilde{\Delta}$, the
lift of the triangulation
$\Delta$ to 
  the universal cover ${\widetilde F_g^s}$ of $F_g^s$. The homeomorphism assigns then to each foliation $F$
representing an element of 
  ${\cal BM}(\Delta)$ an affine foliation $F$ whose lift 
   ${\widetilde F}$ to ${\widetilde F_g^s}$  is equipped with the transverse measure obtained by lifting to the face $F_1$
    the transverse measure of $F$ on the image by the covering map, and then scaling the transverse measures on the various faces 
    of $\widetilde{\Delta}$ so that they agree with each other along the boundaries of these faces.

\vskip .3in

\noindent  $\underline{ \bf
\hbox{Appendix B--- $\lambda$-lengths and broken hyperbolic structures.}}$ 
\vskip .2in
\medskip

As illustrated in Propostion~3.4, $\lambda$-lengths give global coordinates on $\widetilde{\cal BH}(\Delta )$.  In this
appendix, we recall and apply some of the machinery of $\lambda$-lengths to broken hyperbolic structures.

\medskip

There are several coordinatizations of spaces of broken and unbroken hyperbolic structures on a punctured surface which
depend upon the Minkowski
inner product $<\cdot ,\cdot >$ on $\real^3$ whose quadratic form is given by
$x^2+y^2-z^2$ in the usual coordinates.   As is well-known, the upper sheet
$${\bf H}=\{ u=(x,y,z)\in\real^3:<u,u> =-1~~{\rm and}~z>0\}$$ of the two-sheeted hyperboloid is isometric to the
hyperbolic plane.  Furthermore, the open positive light cone
$$L^+=\{ u=(x,y,z)\in\real^3:<u,u> =0~~{\rm and}~z>0\}$$ is identified with the collection of all horocycles in ${\bf
H}$ via  the correspondence $u\mapsto h(u)=\{ w\in{\bf H}:<w,u>=-1\}$.
Recall that  the lambda length of a pair of horocycles $\{ h_0,h_1\}$ 
is defined to be the transform $\lambda (h_0,h_1)=\sqrt{2~{\rm exp}~\delta}$.  Taking this particular transform renders the
identification
$h$ geometrically natural in the sense that 
$\lambda (h(u_0),h(u_1))=\sqrt{-<u_0,u_1>}$, for $u_0,u_1\in L^+$ as one can check.

\medskip

Three useful lemmas (with computational proofs, which we do not here reproduce) are as follows:

\medskip 

\noindent {\bf Lemma B.1} [Pe1;Lemma 2.4] {\it Given three linearly independent rays $r_0,r_1,r_2\subseteq L^+$ from the origin and given
three numbers $\lambda _0,\lambda _1,\lambda _2\in\real_+$, there are unique $u_i\in r_i$, for $i=0,1,2$ so that
$\lambda (h(u_i),h(u_j))=\lambda _k$, where $\{ i,j,k\} =\{ 0,1,2\}$.}

\medskip 

\noindent {\bf Lemma B.2} [Pe1;Lemma 2.3] {\it Given two points $u_0,u_1\in L^+$, which do not lie on a common ray through the origin, and
given two numbers $\lambda _0,\lambda _1\in\real_+$, there is a unique point $v$ on either side of the plane through the origin containing
$u_0,u_1$ satisfying $\lambda (h(v),h(u_i))=\lambda _i$, for $i=0,1$.}

\medskip

\noindent {\bf Lemma~B.3} [Pe1;Proposition 2.8] {\it 
Suppose that $\{ u_i\} _0^2\subseteq L^+$ are linearly independent, let $\gamma (u_i,u_j)$ denote the 
geodesic in ${\bf H}$ with ideal vertices given by the centers of $h(u_i)$ and $h(u_j)$, for $i\neq j$, and define
$$-\lambda _i^2= <u_j,u_k>,~~\alpha _i={{\lambda _i}\over{\lambda _j\lambda _k}},~~for~~\{ i,j,k\} =\{ 0,1,2\} .$$
Then $2\alpha _i$ is the hyperbolic length along the horocycle $h(u_i)$ between $\gamma (u_i,u_j)$ and $\gamma (u_i,u_k)$,
for $\{ i,j,k\} =\{ 0,1,2\}$.
}

\medskip

Armed with these lemmas, it is not difficult to give the basic coordinatization of decorated Teichm\"uller space by lambda lengths:

\medskip

\noindent {\bf Theorem B.4} [Pe1;Theorem 3.1] {\it Fix an ideal triangulation $\Delta$ of $F_g^s$.  Then 
the assignment of $\lambda$-lengths $\tilde{\cal T}(F)\to R^\Delta$is a homeomorphism onto.}

\medskip

\noindent\dem
We must describe an inverse to the mapping and thus give the construction of a decorated hyperbolic structure from 
an assignment of putative $\lambda$-lengths.  To this end, consider the topological universal cover $\tilde F_g^s$ of $F_g^s$ and the lift 
$\tilde\Delta$ of $\Delta$ to
$\tilde F_g^s$; to each component arc of $\tilde\Delta$ is associated the lambda length of its projection.    

The proof proceeds by
induction, and for the basis step, choose any face $T_0$ of
$\tilde\Delta$ and any ideal triangle $t_0$ in
${\bf H}$.  The ideal vertices of $t_0$ determine three rays in $L^+$, so by Lemma~A.1, there are three well-defined points in $L^+$
realizing the putative $\lambda$-lengths on the edges of $T_0$.  (In effect, the basis step of normalizing a triangle ``kills'' the
conjugacy by the M\"obius group in the definition of Teichm\"uller space.)  Of course, the triple of points in $L^+$ corresponds by affine
duality to a triple of horocycles, one centered at each ideal vertex of $t_0$, i.e., a ``decoration'' on $t_0$.

To begin the induction step, consider a face $T_1$ adjacent to $T_0$ across an arc in $\tilde\Delta$.  The two ideal
points $u,v$ which $T_0$ and $T_1$ share have been lifted to $L^+$ in the basis step, and we let $w$ denote the third ideal point of $T_0$
and let
$\Pi$ denote the plane through the origin determined by $u,v$.  According to Lemma~A.2, there is a unique lift $z$ of the third ideal point
of
$T_1$ to $L^+$ on the side of
$\Pi$ not containing $w$ which realizes the putative $\lambda$-lengths.  Again, $u,v,z$ gives rise via affine duality to another decorated
triangle $t_1$ in ${\bf H}$ sharing one edge and two horocycles with $t_0$.

One continues in this manner serially applying Lemma~A.2 to produce a collection of decorated triangles in ${\bf H}$, where any two
triangles have disjoint interiors (because of our choice of the side of the plane in Lemma~A.2).  Thus, our construction gives an injection
$\tilde F_g^s\to {\bf H}$, and we next show that in fact this mapping is also a surjection.  To this end, note first that the inductive
construction has an image which is open in ${\bf H}$ by construction.  According to Lemma~A.3, there is some $\varepsilon >0$ so that each
horocyclic arc inside of each triangle has length at least $\varepsilon$; indeed, there are only finitely many values for such lengths
because the surface is comprised of finitely many triangles.  Thus, each application of the inductive step moves a definite amount along each
horocycle, and it follows easily that the induction has an image which is closed as well.  It follows from connectivity that $\tilde
F_g^s\to{\bf H}$ is a homeomorphism, so $\tilde\Delta$ maps to a tesselation of ${\bf H}$, i.e., a locally finite collection of geodesics
decomposing ${\bf H}$ into ideal triangles.

Following Poincar\'e, the hyperbolic symmetry group of this tesselation is the required (normalized) Fuchsian group $\Gamma$ giving a point
of Teichm\"uller space, and the construction likewise provides a decoration on the quotient ${\bf H}/\Gamma$ as required.
\cqfd

\medskip

By a ``sector'' of $\Delta$ we mean a triangle-vertex pair, and we let $\Sigma =\Sigma (\Delta )$ denote the collection of sectors.  To each
sector is assigned half the hyperbolic length, or the ``h-length'', of the corresponding horocyclic segment as given by Lemma~A.3.  Suppose
that two decorated ideal triangles
$t_0,t_1$ meet along a common edge
$e$, and suppose that the h-lengths of the sectors incident on the endpoints of $e$ in $t_0$ are $\alpha ,\beta$ and in $t_1$ are $\gamma
,\delta$; it follows directly from Lemma~A.3 that $\alpha\beta=\gamma\delta$, the so-called ``coupling equation'' on h-lengths, and we have:

\medskip

\noindent {\bf Corollary B.5} [Pe1;Proposition~3.5]~~{\it  Decorated Teichm\"uller space is identified with the quadric variety in
$\real_+^\Sigma$ determined by the coupling equations on h-lengths.}
\cqfd

\medskip

The reader will recognize the logs of the h-lengths as the transform in (4.1.1).

\medskip

In order to apply the foregoing to broken hyperbolic structures, we next recall several definitions from [OP1], whose
ideas originate in [Th].

\medskip

\noindent {\bf Definition B.6 (Stretch map).---} For $i=1,2$, let 
$(S_i,\Delta_i)$ be a pair consisting of a
 hyperbolic surface $S_i$ with punctures,
   equipped with a geodesic triangulation $\Delta_i$ whose
   vertices are at the  
 punctures, and such that each face  of $\Delta_i$ is isometric to an
 ideal triangle. Let $K>1$. A map $f:(S_1,\Delta_1)\to (S_2,\Delta_2)$ is said 
to be
 a {\it stretch map}  with factor $K$ if the restriction of $f$ to each
face of $S_1$ is a $K-$Lipschitz stretch map of the following form :
 
 \medskip

 \noindent {\bf Definition B.7 ($K-$Lipschitz stretch map).---} Let $K>1$ and consider two hyperbolic ideal triangles $T_1$ and
 $T_2$.
  Then, a
  {\it $K-$Lipschitz stretch map from $T_1$ to $T_2$} is
    is a $K-$Lipschitz homeomorphism
   such that on each side of  $T_1$, the map acts as a
 homothety,
 multiplying
 arclength by $K$, and such that on every closed disk 
 contained in the interior of
 the triangle, the global Lipschitz constant of the restricted map is strictly 
less then $K$.
 For $K=1$,  a $K-$Lipschitz stretch map  is an isometry.
 For $K<1$, a 
 $K-$Lipschitz stretch map is defined as the inverse map of a
$1/K-$Lipschitz stretch map of factor $1/K$.

\medskip

Suppose that $\phi : \pi_1(F_g^s)\to \real_+$ is a homomorphism to the multiplicative group of positive reals.
(The following definition is adapted from [OP1], where a similar structure is
defined on the universal abelian cover of $F_g^s$, instead of on the universal 
cover as we do here.)

\medskip

\noindent {\bf Definition B.8 (Covering hyperbolic structure).---}  A 
{\it covering hyperbolic
structure for $(F_g^s,\Delta)$ with stretch homomorphism 
$\phi$ } is a hyperbolic structure $h$ on $\widetilde{F_g^s}$
 such that the lift of the triangulation $\Delta$ to  
 $\widetilde{F_g^s}$ is a geodesic triangulation $\widetilde{\Delta}$, 
 and such that each covering translation $\alpha\in\pi_1(F_g^s)$ of the cover
 $\widetilde{F_g^s}\to F_g^s$ restricts to a stretch map on each triangle complementary to $\tilde\Delta$ whose factor is $\phi(\alpha)$.

\medskip

 \noindent {\bf Construction B.9 (Developing image of broken hyperbolic
 structure).---}
 Let $H$ be a broken hyperbolic structure on $(F_g^s,\Delta)$, and take the corresponding $\lambda$-lengths
in $\real _+^P$, where $P$ is the collection of triangle-edge pairs of $\Delta$.  Choose a face $T_0$ of $\Delta$ to normalize and use the
$\lambda$-lengths to define a lift to $L^+$ as in the basis step of the proof of Theorem~B.4.  Suppose that $T_1$ is another face of $\Delta$
meeting $T_0$ along the edge $e$, so $\lambda (T_1,e)=\sigma (T_0,T_1)\lambda (T_0,e)$, where $\sigma$ denotes the homothety scaling factor
in (1.1.1).  Let $(T_1,a),(T_1,b)$ denote the triangle-edge pairs of $T_1$ other than $(T_1,e)$, and apply Lemma~A.2 as before using the
$\lambda$-lengths $\sigma (T_0,T_1)\lambda (T_1,a)$ and $\sigma (T_0,T_1)\lambda (T_1,b)$ to produce the lift $z\in L^+$ of the vertex of
$T_1$ other than the ideal points of $e$.  Notice that the center of the horocycle $h(z)$ is actually independent of $\sigma (T_0,T_1)$ by
linearity, though the horocycle $h(z)$ itself does depend upon the $\sigma (T_0,T_1)$.  Continue in this manner as in
Theorem~B.4, to produce a mapping $\tilde F_g^s\to {\bf H}$, whose image is called the {\it developing image} of the
broken hyperbolic structure $H$.

 \medskip

 \noindent {\bf Proposition B.10---} {\it The hyperbolic structure on 
 $\widetilde{F_g^s}$ which is produced by Construction B.9 is a covering 
hyperbolic
 structure.}
\cqfd 
 
 \medskip

\noindent{\bf Theorem B.11} {\it A broken hyperbolic structure on 
a punctured surface of finite
topological type gives rise to a complete
hyperbolic structure on the plane; in other words, 
the developing image of this structure is the entire
hyperbolic plane. The universal cover is thus isometric to the hyperbolic
plane with the covering transformations acting by
quasi-isometries.}

\medskip

\noindent\dem
It remains only to prove that the developing image of a fixed broken hyperbolic structure is the entire hyperbolic plane.  Consider a
triangle
$t_1\subset{\bf H}$ and its translate
$t_{n+1}$ under a parabolic fixing a vertex $u$ of $t_1$, and let $t_1,t_2,\ldots ,t_{n+1}$ be the consecutive faces of
$\Delta$ in the developing image from $t_1$ to $t_{n+1}$.  Since
 the holonomy homomorphism
of a broken hyperbolic structures is trivial around the punctures,
 we have $\prod _{i=1}^{n} \sigma
(t_i,t_{i+1})=1$, so in the developing image the horocyclic segment in $t_1$ centered at $u$ has the same hyperbolic length as the horocyclic
segment in $t_{n+1}$ centered at $u$ by Lemma~A.3.  Now, let $N$ majorize the number of intersections of a horocyclic leaf of a horocyclic
foliation with
$\Delta$ (so $N$ can be estimated merely in terms of the topology of the surface), and let $\varepsilon$ be the smallest h-length of the
broken hyperbolic structure in any sector of $\Delta$.  Thus, traversing $N$ consecutive faces sharing an ideal point in the developing image
moves at least a distance $\varepsilon$ along the horocycle.  It follows that the developing image is the entire hyperbolic plane.
\cqfd

\medskip

This is an interesting class of geometric structures generalizing hyperbolic structures on surfaces: the universal cover admits the metric
structure of the hyperbolic plane with the covering transformations acting by quasi-isometry for this metric.  In our investigations here,
we have restricted attention to such structures relative to a specified ideal triangulation; we wonder if there is an application of the
generalized convex hull construction in Section~7  of [Pe3] (which applies to discrete and radially dense subsets of
$L^+$) giving a corresponding cell decomposition of an appropriate space of these generalized hyperbolic structures in
the spirit of Section 5 of [Pe1]. 

\medskip 

Finally, given a tesselation $\tau$ of ${\bf H}$, define a {\it broken hyperbolic structure} of
$({\bf H},\tau )$ to be a Riemannian metric of constant curvature -1 on ${\bf H}-\tau$, such that each face is isometric to an ideal
triangle, where we furthermore require that metrics on adjacent faces $t_0,t_1$ are related by a homothety scaling factor $\sigma
(t_0,t_1)$, where there is no completeness condition and for the no horocyclic holonomy condition, we require that there
is some
$K>1$ so that for each edge
$e$ of $\tau$ there is some natural number $n=n_e$ so that if $e=e_1,e_2,\ldots, e_{n+1}$ are consecutive edges of
$\tau$ sharing an ideal point with consecutive faces $f_1,f_2,\ldots ,f_n$, where $e_i,e_{i+1}$ lie in the frontier of
$f_i$, for each $i=1,\ldots ,n$, then
$K^{-1}~~<~~\prod _{i=1}^n \sigma (f_i,f_{i+1})~~<~~K$.  As in Theorem~B.4 or Proposition~3.4, $\lambda$-lengths on triangle-edge pairs of
$\tau$ give global coordinates on the space of broken hyperbolic structures on ${\bf H}$.  Furthermore as in
Theorem~B.11, the developing image of
${\bf H}$ is the entire hyperbolic plane.
  This gives a new class of homeomorphisms of the circle as in
Section 6 of [Pe3].

\eject
\centerline{References}

\medskip

\medskip

\noindent [HO] A. Hatcher and U. Oertel, Affine lamination 
spaces for surfaces, Pacific J.
of Math., 154  (1989), no. 1, 41-73.

\medskip

\noindent [OP1] U. Oertel and A. Papadopoulos, Affine foliations and covering hyperbolic structures,
 Manuscripta Math. 104 (2001), 383-406.

\medskip

\noindent  [OP2] U. Oertel and A. Papadopoulos, 
Broken hyperbolic structures and
  affine foliations on surfaces, preprint 2003.

\medskip

\noindent  [Pa]   A. Papadopoulos, Piecewise-linear coordinates for affine foliations on surfaces. Milan J. Math. 70 (2002), 265--290.

\medskip

\noindent  [PP]   A. Papadopoulos and R. C. Penner,
On the Weil-Petersson symplectic structure at Thurston's boundary, Trans. AMS 
  335   (1993) no. 2, 891-904.

\medskip

\noindent  [Pe1] R. C. Penner, The decorated Teichm\"uller space of punctured
surfaces, Comm. Math. Phys. 113 (1987), 299-339.
  
  \medskip

\noindent [Pe2] ---, Weil-Petersson volumes. J. Differential Geom. 35 (1992), no. 3, 559--608. 
   \medskip

\noindent [Pe3] ---, Universal constructions in Teichm\"uller theory. Adv. Math 98 (1993), 143--215.
\medskip

\noindent [Th] W. P. Thurston, Some stretch maps between hyperbolic surfaces,
preprint, 1984.

\bigskip

 {Athanase Papadopoulos}
 
 {Institut de Math\'ematiques,}

 {CNRS and Universit\'e Louis Pasteur,}

 {7 rue Ren\'e Descartes,}  
 
 {67084 Strasbourg Cedex, France}

\bigskip

 {R. C. Penner}
 
 {University of Southern California,}

 {Departments of Mathematics and Physics,}

 {Los Angeles, CA 90089, USA}

\end{document}